\documentclass[a4paper,12pt,reqno]{amsart}
\usepackage{amsfonts}
\usepackage{amsmath}
\usepackage{amssymb}
\usepackage[a4paper]{geometry}
\usepackage{mathrsfs}
\usepackage{csquotes}
\usepackage{mathrsfs}
\usepackage[toc,page]{appendix}
\usepackage[colorlinks]{hyperref}
\renewcommand\eqref[1]{(\ref{#1})} %Need with hyperref
\numberwithin{equation}{section}
\theoremstyle{plain}
\newtheorem{thm}{Theorem}[section]
\newtheorem{prop}[thm]{Proposition}

 \newtheorem{exa}[thm]{Example}
\theoremstyle{definition}

\newtheorem{rem}[thm]{Remark}

\newcommand{\G}{\mathbb G}

\def\e[#1]{{\textrm{e}}^{#1}}

\def\G{{\mathbb G}}

\newcounter{tmp}

\begin{document}

   \title[Examples of $L^p$-$L^q$ multipliers of non-differential operators]
 %  {Critical cases of Sobolev inequalities on graded groups and applications to fractional order nonlinear subelliptic equations}
 {A note on spectral multipliers on Engel and Cartan groups}

\author[Marianna Chatzakou]{Marianna Chatzakou}
\address{
  Marianna Chatzakou:
  \endgraf
  Department of Mathematics: Analysis, Logic and Discrete Mathematics
  \endgraf
  Ghent University
  \endgraf
  Krijgslaan 281, Building S8
  \endgraf
  B 9000 Ghent, Belgium
  \endgraf
  {\it E-mail address} {\rm Marianna.Chatzakou@UGent.be}
  }

 \thanks{The author was supported by the FWO Odysseus 1 grant G.0H94.18N: Analysis and Partial Differential Equations. }


\dedicatory{Dedicated to Prof Jacques Dixmier for his 97th birthday on 24 May 1924}

     \begin{abstract}
     The aim of this short note is to give examples of $L^p$-$L^q$ bounded spectral multipliers for operators involving left-invariant vector fields and their inverses, in the settings of Engel and Cartan groups. The interest in such examples lies in the fact that a group does not have to have flat co-adjoint orbits, and that the considered operator is not related to the usual sub-Laplacian. The discussed examples show how one can still obtain $L^p$-$L^q$ estimates for similar operators in such settings. As immediate consequences, one gets the corresponding Sobolev-type inequalities and heat kernel estimates. 
     \end{abstract}
     \maketitle

 %      \tableofcontents

\section{Introduction}
\label{SEC:intro}

In the literature, one can find many results on the boundedness of spectral multipliers in different settings. However, in most investigations one considers functions of operators having nice spectral or geometric properties, such as functions of the Laplacian or the sub-Laplacian in the Riemannian and sub-Riemannian settings. 

The purpose of this short note is to show that one can effectively work with other examples. For this, we are interested in the analysis on nilpotent Lie groups, and in the $L^p$-$L^q$ estimates for spectral multipliers there. 

In general, in the setting of stratified or graded (nilpotent) Lie groups, one works with the so-called Rockland operators, which are the left-invariant homogeneous hypoelliptic partial differential operators.
However, here we are interested in looking at spectral multipliers of other, more peculiar operators.

Thus, we will show the $L^p$-$L^q$ boundedness of spectral multipliers $\phi(D)$ for $D$ being a non-Rockland operator in the setting of Engel and Cartan groups. These groups, besides the Heisenberg group, are examples of stratified Lie groups which are well-understood. At the same time, they provide for geometric and algebraic properties rather different from those of the Heisenberg group, e.g. {\em non flat co-adjoint orbits, step higher than 2.}

The results that we show are an application of the groups' analogue of the H\"{o}r\-mander $L^p$-$L^q$ multiplier theorem  \cite[Theorem 1.11]{Hor60}. As a straightforward application, we then recover Sobolev-type embedding inequalities and heat kernel estimates. 
Similar results have been shown for the cases where $D$ is the sub-Laplacian on compact Lie groups, or on the Heisenberg group $\mathbb{H}^n$ in \cite{AR20}, and later on, for Rockland operators $D$ on any graded group in \cite{RR18}. 

   The following theorem (\cite[Theorem 1.2]{AR20}) is a generalisation to locally compact groups of the $L^p$-$L^q$ Fourier multiplier theorem proved by H\"{o}rmander (\cite[Theorem 1.11]{Hor60}). In particular, it shows that the $L^p$-$L^q$ norms of spectral multipliers $\phi \left( |\mathcal{L}| \right)$ depend (only) on the growth rate of the spectral projections of the  operator $\mathcal{L}$.
Recall that by a \textit{left Fourier multiplier} in this setting we mean an operator $\mathcal{L}$ acting on the Schwartz space $\mathcal{S}(\G)$ via
\begin{equation}\label{l.f.m.symbol}
\mathcal{F}_{\G}(\mathcal{L}f)(\pi)=\sigma_{\mathcal{L}}(\pi)\widehat{f}(\pi)\,,\quad \pi \in \widehat{{\G}}\,,
\end{equation}
for an operator $\sigma_{\mathcal{L}}(\pi)$ acting on the representation space of $\pi$. 

\begingroup
\setcounter{tmp}{\value{thm}}% store current value of theorem counter
\setcounter{thm}{0} %assign desired value to theorem counter
\renewcommand\thethm{\Alph{thm}}% locally redefine the representation of the theorem counter
\begin{thm} 
	\label{Theorem 1.1,AR}
	Let $\G$ be a locally compact separable unimodular group and let $\mathcal{L}$ be a left Fourier multiplier on $\G$. Assume that $\phi$ is monotonically decreasing continuous function on $[0, \infty)$ such that 
	\[
	\phi(0)=1\,,
	\]
	\[
	\lim_{s \rightarrow \infty} \phi(s)=0\,.
	\]
	Then we have the inequality 
	\begin{equation*}
	\label{Theorem 1.1,ineqiality}
	\| \phi \left( | \mathcal{L}| \right)\|_{L^p(\G) \rightarrow L^q(\G)} \lesssim \sup_{s>0} \phi(s) \left[ \tau \left(E_{(0,s)}\left(|\mathcal{L}| \right) \right)\right]^{\frac{1}{p}-\frac{1}{q}}\,, \quad 1<p \leq 2 \leq q < \infty\,.
	\end{equation*}
\end{thm}
\endgroup

\setcounter{thm}{\thetmp} % restore value of theorem counter

We denote by $E_{(0,s)}\left(|\mathcal{L}| \right)$, the spectral projections associated to the operator $\mathcal{L}$ to the interval $(0,s)$, and $\tau$ is the \textit{canonical trace} on the right group von Neumann algebra $VN_{R}(\G)$. Since $\mathcal{L}$ is affiliated with $VN_{R}(\G)$, its decomposition (see \cite[Theorem 1 on page 187]{Dix81}), allows for a decomposition of the spectral projections  $E_{(0,s)}\left(|\mathcal{L}| \right)$ as
\[
E_{(0,s)}(|\mathcal{L}|)=\bigoplus_{\widehat{{\G}}} E_{(0,s)}\sigma_{\mathcal{L}}(\pi)\,d\mu(\pi)\,,
\] 
where $d\mu(\pi)$ is the Plancherel measure on $\widehat{\G}$, and where $\pi \in \widehat{{\G}}$. As a consequence, see  \cite[Theorem 1 on page 225]{Dix81} the trace $\tau(E_{(0,s)}(\mathcal{L}))$ can be decomposed via 
\begin{equation}
\label{decomp.trace.Dix}
\tau(E_{(0,s)}(\mathcal{L}))=\int_{\widehat{\G}}\tau(E_{(0,s)}(\sigma_{\mathcal{L}}))\,d\mu(\pi)\,.
\end{equation}

\indent To motivate the examples of this paper, let us first recall the corresponding results in the compact case and in the case of the Heisenberg groups $\mathbb{H}^n$. Let $-\Delta_{sub}$ be the positive sub-Laplacian on a compact Lie group $\G$, or in the nilpotent case, on the Heisenberg group $\mathbb{H}^n$. Then, the trace of the spectral projections has the asymptotics
\begin{equation}
\label{compact,heisen.}
\tau \left(E_{(0,s)}\left(-\Delta_{sub} \right) \right)\lesssim s^{Q/2}\,,\quad \textrm{as}\, \quad s \rightarrow \infty\,,
\end{equation}
where $Q$ is the homogeneous dimension in the case of $\mathbb{H}^{n}$, and the Hausdorff dimension generated by the control distance of $\Delta_{sub}$ in the case of the compact Lie group $\G$, respectively. This result has been shown in the compact case in \cite{HK16}, and in the Heisenberg setting in \cite{AR20}. Consequently, an application of Theorem \ref{Theorem 1.1,AR} yields for both cases, for 
$\phi$ as in Theorem \ref{Theorem 1.1,AR}, and $1<p \leq 2 \leq q< \infty$, the inequality
	$$\| \phi \left( -\Delta_{sub} \right)\|_{L^p(\G) \rightarrow L^q(\G)} \lesssim \sup_{s>0} \phi(s)\cdot s^{Q/2}.$$
Recently, in \cite[Theorem 8.2]{RR18}, the authors, employing different machinery from the one we use in this work, proved that \eqref{compact,heisen.} holds true for any graded Lie group of homogeneous dimension $Q$, and, more generally they showed that 
\[
\tau \left(E_{(0,s)}\left(\mathcal{R} \right) \right)\lesssim s^{Q/\nu}\,, \quad \textrm{as}\, \quad s \rightarrow \infty\,,
\]
where $\mathcal{R}$ is a positive Rockland operator of homogeneous degree $\nu$. As a result, they recovered the Sobolev embedding estimates for the case where the Sobolev norm is associated with the sub-Laplacian, or more generally  with any Rockland operator, as in  \cite[Section 4]{Fol75}, and in \cite{FR17} or \cite[Section 4.4.6]{FR16}, respectively.

An interesting question is what happens when we take spectral multipliers of non-Rockland operators, in which case we do not have an easy access to the properties of their spectral projections. So in this paper, we give examples of such applications to the $L^p$-$L^q$ boundedness of the spectral multipliers of some non-Rockland operators on the Engel and Cartan groups. As another application, we get Sobolev type estimates, and time asymptotics for the $L^p-L^q$ norms of the heat kernels of the corresponding heat equations.

  Besides the big amount of work devoted to the case of the Heisenberg group which is of step 2, the same motivating aspects often hold for other graded Lie groups. In particular, we look at two higher step groups here, the Engel and the Cartan group, whose
  representations have been explicitly studied by Dixmier \cite{Dix57}, and whose heat kernels have been also analysed in explicit detail in \cite{BGR10}. We show how this knowledge, combined with Theorem \ref{Theorem 1.1,AR}, implies further regularity estimates.
  
 Let us note that in \cite[Example 7.5]{AR20} the proof of the analogous results to those in this note, in the setting of the Heisenberg group $\mathbb{H}^n$, relies on the fact that the (global) symbol of the canonical sub-Laplacian on $\mathbb{H}^n$, is the harmonic oscillator, and thus the eigenfunctions (the well-known Hermite functions), as well as the corresponding eigenvalues are known explicitly. However, in the  examples of this paper we proceed without knowing the exact eigenvalues of the symbol of the Fourier multiplier in each setting using  estimates for the asymptotic behaviour of the eigenvalue counting function $N(s)$ for anharmonic oscillators $-\frac{d^2}{du^2}+V(u)$ on $L^2(\mathbb{R})$, using some analysis of \cite{CDR18}.
 
We refer to \cite{AK19} for physical analysis of some step 3 groups.
The details of the pseudo-differential theory on Engel and Cartan groups in a rather explicit form have appeared in \cite{Cha20,Cha21}.
 Representations theory of nilpotent Lie groups can be found e.g. in \cite{CG90}, but for Engel and Cartan groups we can also rely on explicit calculation by Dixmier \cite{Dix57}.
 The operators considered in this paper are not covered by more general pseudo-differential theorems, such as Mihlin-H\"ormander theorems for Fourier multipliers \cite{FR21} or for general pseudo-differential operators \cite{CDR21} on graded Lie groups, see also its analogue on compact Lie groups \cite{DR19}.
 As usual, operators on groups like Engel and Cartan group have general importance in view of the Rothschild-Stein lifting theorem \cite{RS77}, and they also fall within the scope of the general analysis on stratified \cite{Fol75} and homogeneous groups \cite{FS82}.

\section{The Engel group}
\label{Engel}

   Let $\G$ be a nilpotent Lie group, and let $\widehat{\G}$ be its  unitary dual, that is, the set of all equivalence classes of irreducible, strongly continuous and unitary representations of $\G$. The \textit{Fourier transform} of $\kappa\in L^{1}(\G)$ at $\pi \in \widehat{\G}$ is defined by
  \begin{equation}\label{group,fourier,gen}
\pi(\kappa)= \widehat{\kappa}(\pi):= \int_{\G}\kappa(x) \pi(x)^{*} dx\,,
\end{equation}
where $dx$ is the Haar measure on $\G$, and $\pi^{*}(x)$ the adjoint operator of $\pi(x)$. This defines a linear mapping $\widehat{\kappa}(\pi):\mathcal{H}_{\pi} \rightarrow\mathcal{H}_{\pi}$ on the representation space $\mathcal{H}_{\pi}$ of $\pi$.
   The infinitesimal representation (or the group Fourier transform) of a vector $X \in \mathfrak{g}$ in the Lie algebra of the group at $\pi \in \widehat{\G}$ is the operator $\pi(X)$ on $\mathcal{H}_{\pi}^{\infty}\subset \mathcal{H}_{\pi}$, the subspace of smooth vectors, given via
   \begin{equation}\label{defn.infin}
   \pi(X)v=\partial_{t=0}(\pi(\textnormal{exp}_{\G}(tX))v\,,
   \end{equation}
   where $\textnormal{exp}_{\G}:\mathfrak{g}\rightarrow \G$ is the exponential mapping from $\mathfrak{g}$ to $\G$. By setting $\pi(X^{\alpha})=\pi(X)^{\alpha}$ we extend the group Fourier transform to the universal enveloping algebra $\mathfrak{U}(\mathfrak{g})$.

We now turn to the specifics of the Engel group.
For more details of explicit elements of the pseudo-differential theory on Engel and Cartan groups we can refer to \cite{Cha20,Cha21}.

The Engel Lie algebra $\mathfrak{l}_4=\text{span}\{I_1,I_2,I_3,I_4\}$ is the 3-step nilpotent Lie algebra, with the following non-trivial commutator relations:
\[
[I_1,I_2]=I_3\,, \quad [I_1,I_3]=I_4\,.
\]
The Lie algebra $\mathfrak{l}_4$ can be decomposed as
\begin{equation}\label{decomp,eng}
\mathfrak{l}_4=V_1\oplus V_2 \oplus V_3,
\end{equation}
with 
\[
V_1=\text{span}\{I_1,I_2\}\,, V_2=\text{span}\{I_3\}\,, \textrm{and}\,, V_3=\text{span}\{I_4\}\,,
\]
such that $[V_i,V_j] \subset V_{i+j}$.
It is a stratified Lie algebra, with $V_1$ generating the whole of  $\mathfrak{l}_4$. The corresponding Lie group is called   the \textit{Engel group}, denoted by $\mathcal{B}_4$. It is 
a homogeneous Lie group, and the natural dilations on its Lie algebra are given by 
\begin{equation}\label{Dil,eng}
D_r(I_1)=rI_1\,, D_r(I_2)=rI_2\,,D_r(I_3)=r^2 I_3\,, \textrm{and}\, D_r(I_4)=r^3 I_4\,, \quad r>0\,.
\end{equation}
The group $\mathcal{B}_4$ can be identified with the manifold $\mathbb{R}^4$ endowed with the group law:
\begin{eqnarray*}
\lefteqn{(x_1,x_2,x_3,x_4) \times (y_1,y_2,y_3,y_4)}\\
&:=& (x_1+y_1,x_2+y_2,x_3+y_3-x_1y_2,x_4+y_4+\frac{1}{2}x_1^2y_2-x_1y_3)\,.
\end{eqnarray*}
A canonical basis of the Lie algebra $\mathfrak{l}^4$ was calculated in \cite[Section 3.2]{BGR10}, and can be taken to be
\begin{equation}\label{left.inv.eng}
\begin{split}
X_1(x)&= \frac{\partial}{\partial x_1}\,, \quad X_2(x)=\frac{\partial}{\partial x_2}-x_1 \frac{\partial}{\partial x_3}+ \frac{x^{2}_{1}}{2}\frac{\partial}{\partial x_4}\,, \\
            X_3(x)&= \frac{\partial}{\partial x_3}-x_1 \frac{\partial}{\partial x_4}\,,\quad X_4(x)= \frac{\partial}{\partial x_4}\,,
\end{split}
\end{equation}
for $x=(x_1,x_2,x_3,x_4) \in \mathbb{R}^4$. 
The exponential map $\text{exp}_{\mathcal{B}_4}$ is the identity map, giving the identification
\begin{equation}\label{expB_4}
(x_1,x_2,x_3,x_4)=\text{exp}_{\mathcal{B}_4}(x_1X_1+x_2X_2+x_3X_3+x_4X_4).
\end{equation}
The Haar measure on $\mathcal{B}_4$ is the Lebesgue measure on  $\mathbb{R}^4$. In view of \eqref{expB_4},  the dilations \eqref{Dil,eng} yield the group dilations,  
\[
D_r(x_1,x_2,x_3,x_4)=(r x_1,r x_2, r^2 x_3, r^3 x_4)\,, \quad r>0\,,
\]
and the homogeneous dimension of $\mathcal{B}_4$ is $Q_{\mathcal{B}_4}=1+1+2+3=7$.
The left-invariant sub-Laplacian on the Engel group $\mathcal{B}_4$ is homogeneous of order 2, and given by 
$$
\mathcal{L}_{\mathcal{B}_4}= X_{1}^{2}+X_{2}^{2}
=\frac{\partial^2}{\partial x_{1}^{2}}+ \left(\frac{\partial}{\partial x_2}-x_1 \frac{\partial}{\partial x_3}+ \frac{x^{2}_{1}}{2}\frac{\partial}{\partial x_4} \right)^2.
$$

The representations of the Engel group are known from Dixmier \cite[p.333]{Dix57}.
\begin{prop}\label{Dixmier, Engel}
	The unitary dual of $\mathcal{B}_4$ is $\widehat{\mathcal{B}_4}= \{\pi_{\lambda, \mu}\arrowvert \lambda\neq 0, \mu,\lambda \in \mathbb{R}\}$. For each $(x_1,x_2,x_3,x_4) \in \mathcal{B}_4\,,\pi_{\lambda, \mu}(x_1,x_2,x_3,x_4)$ acts on $L^2(\mathbb{R},\mathbb{C})$ via
\[
\pi_{\lambda, \mu}(x_1,x_2,x_3,x_4)h(u)\equiv \exp \left(i \left(-\frac{\mu}{2\lambda}x_2+\lambda x_4-\lambda x_3 u + \frac{\lambda}{2}x_2 u^2\right)\right)h(u+x_1)\,.
\]
\end{prop}

By \cite[Subsection 3.2.3]{BGR10}, the infinitesimal representations (the symbols) of the basis of $\mathfrak{l}_4$ are the operators acting on $\mathcal{S}(\mathbb{R})$ given by
\begin{equation}
\label{inf1,2,eng}
\pi_{\lambda, \mu}(X_1)=\frac{d}{du},\; \pi_{\lambda, \mu}(X_2)=\bigg( -\frac{i\mu}{2\lambda} +\frac{i}{2}\lambda u^2 \bigg ),\;  \pi_{\lambda, \mu}(X_3)=-i\lambda u,\;  \pi_{\lambda, \mu}(X_4)=i \lambda.
\end{equation}
It follows then that the symbol of the sub-Laplacian $\mathcal{L}_{\mathcal{B}_4}$ is well-defined as an operator acting on 
the Schwartz space $\mathcal{S}(\mathbb{R})$, and is given by
$$
	\pi_{\lambda, \mu}(\mathcal{L}_{\mathcal{B}_4})=\pi_{\lambda, \mu}(X_1^2+X_2^2)= \left(\frac{d}{du} \right)^2+ \left( -\frac{i\mu}{2\lambda}+ \frac{i \lambda u^2}{2}\right)^2
	= \frac{d^2}{du^2} -\frac{1}{4} \left(\lambda u^2 -\frac{\mu}{\lambda} \right)^2\,.
$$

We now show how Theorem \ref{Theorem 1.1,AR} can be applied to operators which are not necessarily differential.

\begin{exa}[$L^p$-$L^q$ multipliers on the  Engel group]
	\label{multipliers,engel}
	Assume that $\phi$ is monotonically decreasing continuous function on $[0, \infty)$ such that 
	\[
	\phi(0)=1\,,
	\]
	\[
	\lim_{s \rightarrow \infty} \phi(s)=0\,.
	\]
	Then, for
	\[
	A=- \left(X_1^2+X_2^2+X_3^2+X_4^2+X_4^{-2} \right)\,,
	\]
	with $X_i$, $i=1, \ldots,4$, being the left-invariant vector fields on $\mathcal{B}_4$, given in \eqref{left.inv.eng}, we have the inequality 
	\[
	\| \phi(A) \|_{L^{p}(\mathcal{B}_4) \rightarrow L^{q}(\mathcal{B}_4)} \lesssim \sup_{s>0}\phi(s)\cdot s^{\frac{3}{r}}\,,
	\]
	where $\frac{1}{r}=\frac{1}{p}-\frac{1}{q}$, for $1<p \leq 2 \leq q < \infty$. This has to be understood that if the right hand side is finite, so is also the left hand side.
	In addition, for $f \in \mathcal{S}(\mathcal{B}_4)$, the Sobolev-type inequalities
	\begin{equation*}
	\label{for_phi,engel}
	\|(I+A)^{b} f\|_{L^{q}(\mathcal{B}_4)} \leq \| (I+A)^{a}f \|_{L^{p}(\mathcal{B}_4)}
	\end{equation*}
	hold true for $a, b \in\mathbb{R}$, such that $a-b\geq 	\frac{3}{r}$, and for $1<p \leq 2 \leq q < \infty$.
\end{exa}
Before turning on to prove the above let us first note that the Fourier multiplier $X_{4}^{-2}$, if $\lambda \neq 0$, can be realised through its (global) symbol, i.e.,
\[
\mathcal{F}_{\mathcal{B}_4}(X_{4}^{-2}f)(\pi_{\lambda,\mu})=-\lambda^{-2}\widehat{f}(\pi_{\lambda,\mu})\,.
\]
Indeed, we can write for $f \in \mathcal{S}(\mathbb{R})$, \[\pi_{\lambda,\mu}(f)=\pi_{\lambda,\mu}(X_{4}^{2}X_{4}^{-2}f)=\pi_{\lambda,\mu}(X_{4}^{2})\pi_{\lambda,\mu}(X_{4}^{-2}f)=(i \lambda)^{2}\pi_{\lambda,\mu}(X_{4}^{-2}f)\,.\]

\begin{proof}[Proof of Example \ref{multipliers,engel}]
	For $A$ as in the statement, using \eqref{l.f.m.symbol} and \eqref{decomp.trace.Dix}, we get
	\begin{equation}
	\label{Dix,trace,eng}
	\tau \left( E_{(0,s)}(A)\right)\approx\int_{\widehat{\mathcal{B}_4}} \tau \left(E_{(0,s)}\left[ \pi_{\lambda,\mu}(A) \right]  \right)    d\lambda\, d\mu\,.
	\end{equation}
	For $\mu \in \mathbb{R}\,, \lambda \neq 0$, the operator $\pi_{\lambda, \mu}(A): \mathcal{S}(\mathbb{R}) \rightarrow \mathcal{S}(\mathbb{R})$, is the (global) symbol of $A$, given by  
	\begin{equation*}\label{R,eng}
	\pi_{\lambda, \mu}(A)=-\frac{d^2}{du^2}+ \frac{1}{4}\left( \lambda u^2- \frac{\mu}{\lambda} \right)^{2}+(\lambda u)^{2}+\lambda^{2}+\lambda^{-2}\,,\quad u \in \mathbb{R}\,,
	\end{equation*}
	as it follows from \eqref{inf1,2,eng}. The operator $\pi_{\lambda, \mu}(A)$ has discrete spectrum (see e.g. \cite{CDR18}), and if $N_{A}(s)$ denotes the number of its eigenvalues not exceeding $s$, then 
	\begin{equation}\label{type1,eng}
	\tau \left(E_{(0,s)}\left[ \pi_{\lambda,\mu}(A) \right] \right) = N_{A}(s)\,.
	\end{equation}
	 Let $\sigma_{A}$ be the Weyl symbol of $A$ (after rescaling). Then, for $\mu \in \mathbb{R}$ and for $\lambda>0$, the latter is given by
	\[
	\sigma_{A}(u,\xi)= \xi^{2}+\frac{1}{4}\left( \lambda u^2- \frac{\mu}{\lambda} \right)^{2}+(\lambda u)^{2}+\lambda^{2}+\lambda^{-2}\,.
	\]
	Then,  \cite[Theorem 3.2]{BBR96} implies that 
	\begin{equation}
	\label{BBR, engel}
	N_{A}(s) \lesssim C \int_{ \sigma_{A}(u,\xi)<s} 1 \cdot du\,d\xi\,, \quad \textrm{as}\quad s\rightarrow \infty\,.
	\end{equation}
	Assuming that $\sigma_{A}(u,\xi)<s$, it follows immediately that
	\begin{equation}
	\label{for_lamda,xi,engel}
	s^{-1/2}< \lambda<s^{1/2}\,,\quad |\xi|<s^{1/2}\,.
	\end{equation}
	and
	\begin{equation}
	\label{sigma<s,engel}
	\frac{1}{4}\left( \lambda u^2-\frac{\mu}{\lambda}\right)^2, (\lambda u)^2<s\,.
	\end{equation}
	Now, by \eqref{for_lamda,xi,engel} and \eqref{sigma<s,engel}, we get 
	$ |u|< \frac{s^{1/2}}{\lambda}< \frac{s^{1/2}}{s^{-1/2}}=s$,
	and
	\begin{equation*}
	-2s^{1/2} \cdot s^{1/2}<-2\lambda s^{1/2}<2(\lambda^2 u^2- \lambda s^{1/2})<\mu< 2(\lambda^{2}u^{2}+\lambda s^{1/2})<2(s+s^{1/2}\cdot s^{1/2})\,.
	\end{equation*}
	In particular, we have $-2s<\mu<4s$. Then, by \eqref{BBR, engel}, and since $|u|<s, |\xi|<s^{1/2}$, we have $N_{A}(s) \lesssim C \cdot s^{3/2}$, for $s \rightarrow \infty$, whereas by \eqref{type1,eng}, the integration on the dual $\widehat{\mathcal{B}_4}$ in \eqref{Dix,trace,eng} for $\lambda \in (s^{-1/2},s^{1/2})$, $\mu \in (-2s,4s)$, becomes
	\begin{align*}
	\tau \left( E_{(0,s)}(A)\right)&\approx\int_{s^{-1/2}}^{s^{1/2}} \int_{-2s}^{4s} \tau \left(E_{(0,s)}\left[ \pi_{\lambda,\mu}(A) \right]  \right)    d\lambda\, d\mu\\
	& \lesssim C \cdot \int_{s^{-1/2}}^{s^{1/2}} \int_{-2s}^{4s} s^{3/2} d\lambda\,d\mu \lesssim C^{'}\cdot s^3\,,\quad s \rightarrow \infty\,.
	\end{align*}
	Summarising we have $ \tau \left( E_{(0,s)}(A)\right) \lesssim C^{'} \cdot s^3$, for $s \rightarrow \infty$. Then, for $\phi$ as in the statement, an application of Theorem \ref{Theorem 1.1,AR}, yields
	\[
	\|\phi\left( A \right)\|_{L^{p}(\mathcal{B}_4) \rightarrow L^{q}(\mathcal{B}_4)} \lesssim C^{'} \cdot \sup_{s>0} \phi(s) \cdot  \tau \left( E_{(0,s)}(A)\right)^{\frac{1}{p}-\frac{1}{q}} \lesssim C^{'} \cdot \sup_{s>0}\phi(s) \cdot  s^{\frac{3}{r}}\,,
	\]
	where $\frac{1}{r}=\frac{1}{p}-\frac{1}{q}$, and $1<p \leq 2 \leq q < \infty$, and this shows the first claim in Example \ref{multipliers,engel}.  Now, if we choose 
	\begin{equation}
	\label{def_of_phi,eng}   
	\phi=\phi(s):=\frac{1}{(1+s)^{a-b}}\,,\quad s>0\,,
	\end{equation}
where $a-b \geq \frac{3}{r}$, then for $f \in \mathcal{S}(\mathcal{B}_4)$,
	\[
	\|(I+A)^{-a} (I+A)^{b} f\|_{L^{q}(\mathcal{B}_4)} \lesssim \|f \|_{L^{p}(\mathcal{B}_4)}\,,
	\]
	or, as claimed in the statement,
	\[
	\|(I+A)^{b} f\|_{L^{q}(\mathcal{B}_4)} \lesssim \| (I+A)^{a}f \|_{L^{p}(\mathcal{B}_4)}\,.
	\]
	The proof is complete.
\end{proof}
We note that, in the above example, the operator $A$ can be generalised as being any non-Rockland operator of the form
\[
A=-(X_{1}^{2}+ X_{2}^{2n_2}\pm X_{3}^{2n_3}\pm X_{4}^{2n_4}\pm X_{4}^{-2n_5})\,,
\]
where $n_i \geq 2$, $n_i \in \mathbb{N}$ for $i=2,\ldots,5$, if we choose the positive sign of the corresponding left-invariant operator if $n_i$ is of the form $n_i=1+2k_i$, $k_i \geq 1$, and the negative sign, otherwise. 
\begin{rem}(Sobolev-type estimates on the Engel group)
	\label{Sob,emb2,eng}
	For $b=0$ in the definition of the function $\phi$ in \eqref{def_of_phi,eng}, we obtain Sobolev-type embedding inequalities of the form
	\[
	\|f\|_{L^{q}(\mathcal{B}_4)} \leq C \| (I+A)^a f\|_{L^{p}(\mathcal{B}_4)}\,,
	\]
	for $a \geq \frac{3}{r}$ and $1<p \leq 2 \leq q < \infty$.
\end{rem}
\begin{rem}(The $A$-heat equation)
	\label{heat,eng}
	For the operator $A$ as in Example \ref{multipliers,engel}, consider the heat equation 
	\begin{equation}
	\label{heat.eq,eng}
	\partial_t u+A u=0\,,\quad u(0)=u_0\,.
	\end{equation}
	It is not difficult to check that, for each $t>0$, the function $u(t)=u(t,x):=e^{-tA}u_0$, is a solution of the initial value problem \eqref{heat.eq,eng}. For $t>0$, also set $\phi=\phi(s):=e^{-ts}$. Then $\phi$ satisfies the assumptions stated before in Theorem \ref{Theorem 1.1,AR}, and as an application of the Example \ref{multipliers,engel} we get
	\begin{equation}\label{norm.op.heat.eng}
	\| e^{-tA}\|_{L^p(\mathcal{B}_4) \rightarrow L^q(\mathcal{B}_4)} \lesssim C \cdot \sup_{s>0}e^{-ts} \cdot s^{\frac{3}{r}}\,,
	\end{equation}
	where $\frac{1}{r}=\frac{1}{p}-\frac{1}{q}$, and $1<p \leq 2 \leq q < \infty$. Using techniques of standard mathematical analysis, we see that 
	\[
	\sup_{s>0} e^{-ts} \cdot s^{\frac{3}{r}}=\left( \frac{3}{tr} \right)^{\frac{3}{r}}e^{-\frac{3}{r}}:=C_{p,q} t^{-3\left(\frac{1}{p}-\frac{1}{q}\right)}\,.
	\]
	Indeed, if we consider the function $g(s)= e^{-ts} \cdot s^{\frac{3}{r}}$, then $g^{'}(s)=s^{\frac{3}{r}-1}e^{-ts} \left( \frac{3}{r}-st \right)$, whereas $g^{'}\left( \frac{3}{rt} \right)=0$, and $g^{''}\left( \frac{3}{rt} \right)<0$. Then, by \eqref{norm.op.heat.eng}, we obtain
	\[
	\|u(t, \cdot) \|_{L^q(\mathcal{B}_4)} \lesssim C_{p,q} t^{-3\left(\frac{1}{p}-\frac{1}{q}\right)} \|u_0\|_{L^{p}(\mathcal{B}_4)},\; 1<p \leq 2 \leq q< \infty\,,
	\]
	where the last equations yields the time decay for the solution of the heat equation in this setting. 
\end{rem}

\section{The Cartan group}
\label{Cartan}
The Cartan Lie algebra $\mathfrak{l}_5=\text{span}\{I_1,I_2,I_3,I_4,I_5\}$ is the 3-step nilpotent Lie algebra, with the following nontrivial commutator relations:
\[
[I_1,I_2]=I_3\,, [I_1,I_3]=I_4\,, [I_2,I_3]=I_5\,.
\]
The Lie algebra $\mathfrak{l}_5$ can be decomposed as 
 \begin{equation}\label{decomp,cart}
\mathfrak{l}_5=V_1\oplus V_2 \oplus V_3,
\end{equation}
with\[
V_1=\text{span}\{I_1,I_2\}\,, V_2=\text{span}\{I_3\}\,, \textrm{and}\,,V_3=\text{span}\{I_4,I_5\}\,,
\]
such that $[V_i,V_j] \subset V_{i+j}$.
It is a stratified Lie algebra, with $V_1$ generating the whole of  $\mathfrak{l}_5$. The corresponding Lie group is called   the \textit{Cartan group}, denoted by $\mathcal{B}_5$. It is 
a homogeneous Lie group, and the natural dilations on its Lie algebra are given by 
\begin{equation*}
D_r(I_1)=rI_1\,, D_r(I_2)=rI_2\,,D_r(I_3)=r^2 I_3\,,D_r(I_4)=r^3I_4\,, \textrm{and}\, D_r(I_4)=r^3 I_5\,, \quad r>0\,.
\end{equation*}
The group $\mathcal{B}_5$ can be identified with the manifold $\mathbb{R}^5$ endowed with the group law
\begin{equation*}
\begin{split}
 &(x_1,x_2,x_3,x_4,x_5) \times (y_1,y_2,y_3,y_4,y_5):=\\
 &  \left(x_1+y_1,x_2+y_2,x_3+y_3-x_1y_2,x_4+y_4+\frac{x_1^2y_2}{2}-x_1y_3, x_5+y_5 + \frac{x_1y_2^2}{2}-x_2y_3+x_1x_2y_2\right).       
\end{split}
\end{equation*}
A canonical basis of $\mathfrak{l}^5$ can be given by (see e.g. \cite[Section 3.3]{BGR10})
\begin{equation}\label{left.inv.car}
\begin{split}
X_1(x)&=\frac{\partial}{\partial x_1}\,, \quad X_2(x)=\frac{\partial}{\partial x_2}-x_1 \frac{\partial}{\partial x_3} 
+ \frac{x_1^2}{2} \frac{\partial}{\partial x_4}+x_1x_2\frac{\partial}{\partial x_5}\,, \\
            X_3(x)&= \frac{\partial}{\partial x_3}-x_1 \frac{\partial}{\partial x_4}-x_2 \frac{\partial}{\partial x_5}\,,\quad  X_4(x)= \frac{\partial}{\partial x_4}\,,\quad X_5(x)=\frac{\partial}{\partial x_5}\,,
\end{split}
\end{equation}
where $x=(x_1,x_2,x_3,x_4,x_5) \in \mathbb{R}^5$.
The Haar measure on $\mathcal{B}_5$ is the Lebesgue measure on  $\mathbb{R}^5$, and the dilations on $\mathcal{B}_5$ are given by
\[
D_r(x_1,x_2,x_3,x_4,x_5)=(rx_1,rx_2,r^2x_3,r^3x_4,r^3x_5)\,.
\] 
The homogeneous dimension of $\mathcal{B}_5$ is $Q_{\mathcal{B}_5}=1+1+2+3+3=10$.
The left-invariant sub-Laplacian on the Cartan group $\mathcal{B}_5$ is homogeneous of order 2, and given by 
$$
\mathcal{L}_{\mathcal{B}_5}= X_{1}^{2}+X_{2}^{2}
=\frac{\partial^2}{\partial x_{1}^{2}}+ \left(\frac{\partial}{\partial x_2}-x_1 \frac{\partial}{\partial x_3} 
+ \frac{x_1^2}{2} \frac{\partial}{\partial x_4}+x_1x_2\frac{\partial}{\partial x_5}\right)^2.
$$
The representations of the Engel group are known from Dixmier \cite[p.338]{Dix57}.
\begin{prop}\label{Dixmier,cartan}
	The unitary dual of $\mathcal{B}_5$ is $\widehat{\mathcal{B}_5}= \{\pi_{\lambda, \mu, \nu}\arrowvert \lambda^2+\mu^2\neq 0, \nu, \lambda, \mu \in \mathbb{R}\}$. For each $(x_1,x_2,x_3,x_4,x_5) \in \mathcal{B}_5$ $\,,\pi_{\lambda, \mu, \nu}(x_1,x_2,x_3,x_4,x_5) $ acts on $L^2(\mathbb{R},\mathbb{C})$ via
\[
\pi_{\lambda, \mu, \nu}(x_1,x_2,x_3,x_4,x_5)h(u)\equiv \exp (i A^{\lambda, \mu, \nu}_{x_1,x_2,x_3,x_4,x_5}(u))h\left( u+\frac{\lambda x_1+\mu x_2}{\lambda^2+\mu^2} \right)\,,
\]
with 
\begin{align*}
A^{\lambda, \mu, \nu}_{x_1,x_2,x_3,x_4,x_5}(u)=& -\frac{1}{2} \frac{\nu}{\lambda^2 +\mu^2}(\mu x_1 -\lambda x_2)+\lambda x_4+ \mu x_5\\
& - \frac{1}{6} \frac{\mu}{\lambda^2 + \mu^2} (\lambda^2 x_1^3+3 \lambda \mu x_1^2 x_2 + 3 \mu^2 x_1 x_2^2 -\lambda \mu x_2^3)\\
& +\mu^2 x_1x_2 u +\lambda \mu (x_1^2-x_2^2) u +\frac{1}{2} (\lambda^2+\mu^2)(\mu x_1-\lambda x_2)u^2\,.
\end{align*}
\end{prop}
The infinitesimal representations (symbols) of the vector fields were calculated in \cite[Subsection 3.3.2]{BGR10}, given by 
\begin{equation}\label{inf2,car}
\begin{split}
\pi_{\lambda, \mu, \nu}(X_1)&=\left(-\frac{i}{2} \frac{\nu \mu}{\lambda^2 + \mu^2} - \frac{i}{2} (\lambda^2 + \mu^2) \mu u^2 + \frac{\lambda}{\lambda^2 + \mu^2} \frac{d}{du} \right)\,,\\
\pi_{\lambda, \mu, \nu}(X_2)&=\left( \frac{i}{2} \frac{ \lambda \nu}{\lambda^2 + \mu^2}+\frac{i}{2} (\lambda^2 + \mu^2) \lambda u^2+\frac{\mu}{\lambda^2 + \mu^2} \frac{d}{du} \right)\,,
\end{split}
\end{equation}
as well as 
\begin{equation}
\label{inf3,4,5,car}
\pi_{\lambda, \mu, \nu}(X_3)=i (\lambda^2 + \mu^2)u\,,\quad \pi_{\lambda, \mu, \nu}(X_4)=i \lambda \\,\quad  \pi_{\lambda, \mu, \nu}(X_5)=i \mu\,.
\end{equation}
Consequently, the symbol of the  sub-Laplacian $\mathcal{L}_{\mathcal{B}_5}$ is given by
	\begin{align*}
	\pi_{\lambda, \mu}(\mathcal{L}_{\mathcal{B}_5})=\pi_{\lambda, \mu}(X_1^2+X_2^2)&=\left(-\frac{i}{2} \frac{\nu \mu}{\lambda^2 + \mu^2} - \frac{i}{2} (\lambda^2 + \mu^2) \mu u^2 + \frac{\lambda}{\lambda^2 + \mu^2} \frac{d}{du} \right)^2\\&+ \left( \frac{i}{2} \frac{\lambda \nu}{\lambda^2 + \mu^2}+\frac{i}{2} (\lambda^2 + \mu^2) \lambda u^2+\frac{\mu}{\lambda^2 + \mu^2} \frac{d}{du} \right)^2\\
	&=\frac{1}{\lambda^2+ \mu^2} \frac{d^2}{du^2}-\frac{1}{4(\lambda^2 + \mu^2)}((\lambda^2 + \mu^2)^2 u^2 + \nu)^2.
	\end{align*}

\begin{exa}[$L^p$-$L^q$ multipliers on the Cartan group]
	\label{multipliers,cartan}
	Assume that $\phi$ is monotonically decreasing continuous function on $[0, \infty)$ such that 
	\[
	\phi(0)=1\,,
	\]
	\[
	\lim_{s \rightarrow \infty} \phi(s)=0\,.
	\]
	Let
	\[
	B=-\left(X_1^2+X_2^2+X_3^2+X_4^2+X_5^2+X_4^{-2}+X_5^{-2} \right),
	\]
	where $X_i$ are the left-invariant vector fields on $\mathcal{B}_5$, given in \eqref{left.inv.car}. Then,
	\[
	\| \phi(B) \|_{L^p(\mathcal{B}_5) \rightarrow L^q(\mathcal{B}_5)} \lesssim C \cdot \sup_{s>0} \phi(s) \cdot s^{\frac{9}{2r}}\,,
	\]
	where $\frac{1}{r}=\frac{1}{p}-\frac{1}{q}$, for $1<p \leq 2 \leq q < \infty$. In addition, for $f \in \mathcal{S}(\mathcal{B}_5)$, we get the Sobolev-type inequalities of the form 
	\begin{equation}
	\label{Sob,cartan}
	\| (I+B)^b f \|_{L^q(\mathcal{B}_5)} \leq \| (I+B)^a f \|_{L^p(\mathcal{B}_5)}\,,
	\end{equation}
	for $a,b \in \mathbb{R}$, such that $a-b\geq \frac{9}{2r}$ and for $1<p \leq 2 \leq q < \infty$.
\end{exa}
As before, the Fourier multiplier $X_{5}^{-2}$, if $\mu\neq 0$, can be realised as
\[
\mathcal{F}_{\mathcal{B}_5}(X_{5}^{-2}f)(\pi_{\lambda, \mu, \nu})=-\mu^{-2}\widehat{f}(\pi_{\lambda, \mu, \nu})\,,
\]
since $\pi_{\lambda, \mu, \nu}(X_5)=i\mu$, and similarly for $X_{4}^{-2}$.
\begin{proof}[Proof of Example \ref{multipliers,cartan}]
	Let $B$ be as in the statement, by using \eqref{l.f.m.symbol} and \eqref{decomp.trace.Dix}, we get
	\begin{equation}
		\label{Dix, cartan}
	\tau \left( E_{(0,s)}(B)\right)\approx\int_{\widehat{\mathcal{B}_5}} \tau \left(E_{(0,s)}\left[ \pi_{\lambda,\mu, \nu}(B) \right]  \right)    d\lambda\, d\mu\,d\nu\,,
	\end{equation} 
	where $d\lambda\,d\mu\,d\nu$ is the Plancherel measure on $\widehat{\mathcal{B}_5}$, see \cite{Dix57}. For $\nu \in \mathbb{R}\,, \lambda,\mu \neq 0$, and $u \in \mathbb{R}$ the operator symbol  $\pi_{\lambda, \mu, \nu}(B) : \mathcal{S}(\mathbb{R}) \rightarrow \mathcal{S}(\mathbb{R})$ of $B$ is  given by 
	\[
	\pi_{\lambda, \mu, \nu}(B)=-\frac{1}{\lambda^2+\mu^2} \frac{d^2}{du^2}+\frac{1}{4}\frac{(\nu+(\lambda^2+\mu^2)^{2}u^2)^2}{\lambda^2+\mu^2}+(\lambda^2+\mu^2)^2 u^2+\lambda^{2}+\mu^{2}+\lambda^{-2}+\mu^{-2}\,,
	\]
	using \eqref{inf2,car} and \eqref{inf3,4,5,car}. Therefore, we have
	\begin{equation}\label{type1,car}
	\tau \left(E_{(0,s)}\left[ \pi_{\lambda,\mu,\nu}(B) \right] \right) = N_{B}(s)\,,
	\end{equation}
	where $N_{B}(s)$ denotes the number of the eigenvalues of $\pi_{\lambda, \mu, \nu}(B)$ that are less than $s>0$.
	Observe that $\pi_{\lambda, \mu, \nu}(B)$ is a rescaled anharmonic oscillator. Now, if we consider the operator $\pi_{\lambda, \mu, \nu}(B^{'})=(\lambda^2+\mu^2) \cdot\pi_{\lambda, \mu, \nu}(B)$,  then by \cite[Theorem 3.2]{BBR96}, for $s \rightarrow \infty$, we have
	\begin{equation}
	\label{BBG,cartan}
	N_{B}(s)=N_{B^{'}}((\lambda^2+\mu^2)\cdot s) \lesssim C \cdot \int_{\sigma_{B^{'}}(u,\xi)< (\lambda^2+\mu^2)\cdot s} 1 \cdot du\,d\xi\,,
	\end{equation}
	where with $N_{B^{'}}((\lambda^2 + \mu^2)\cdot s)$ we have denoted the number of eigenvalues of $\pi_{\lambda, \mu, \nu}(B^{'})$, that are less than $(\lambda^2+\mu^2) \cdot s$, and with $\sigma_{B^{'}}$, the (rescalled) Weyl symbol of the anharmonic oscillator $\pi_{\lambda, \mu, \nu}(B^{'})$ given by
	\[
\sigma_{B^{'}}(u,\xi)=	\xi^2+ (\nu+(\lambda^2+\mu^2)^{2}u^2)^{2}/4+(\lambda^2+\mu^2)^3 u^2 + (\lambda^2+\mu^2)^2+ (\lambda^2+\mu^2) (\lambda^{-2}+\mu^{-2})\,,
	\]
	where we have assumed, without loss of generality, that $\lambda, \mu>0$. Now, if $\sigma_{B^{'}}(u,\xi)<(\lambda^2+\mu^2) \cdot s$, then 
	\begin{equation}
	\label{for_l,m,cartan}
	s^{-1/2}<\lambda, \mu< s^{1/2}\,, \quad \textrm{and}
	\end{equation}
	\begin{equation}
	\label{for_x,xi,nu,cartan}
	(\lambda^2+\mu^2)^{2}u^2<s\,,\quad (\nu+(\lambda^2+\mu^2)^{2}u^2)^{2}<4 s (\lambda^2+\mu^2)\,,\quad \xi^2 < s(\lambda^2+\mu^2)\,.
	\end{equation}
	Now, by \eqref{for_l,m,cartan} and \eqref{for_x,xi,nu,cartan}, we get $|u|<\frac{s^{1/2}}{\lambda^2+\mu^2}<\frac{s^{1/2}}{s^{-1}+s^{-1}}=\frac{s^{3/2}}{2}$, $| \xi |<\sqrt{2} s^{1/2}\cdot s^{1/2}$ and
	\[
	\nu^{2}<(\nu+(\lambda^{2}+\mu^{2})^{2}u^{2})^{2}<4s (\lambda^{2}+\mu^{2})<4s^{2}\,.
	\]
	In particular, we have $| \nu |<2s$. Then, by \eqref{BBG,cartan}, and for $u,\xi$ in the above range, we have 
	\begin{eqnarray}
	\label{N^',cartan}
	N_{B}(s)=N_{B^{'}}((\lambda^2+\mu^2) \cdot s) \lesssim C \cdot s^{3/2} \cdot 2\sqrt{2} s \lesssim 2\sqrt{2} C \cdot s^{5/2}\,,\quad s \rightarrow \infty\,.
	\end{eqnarray}
	Collecting, \eqref{type1,car} and \eqref{N^',cartan}, the integration on the dual $\widehat{\mathcal{B}_5}$ that appears in \eqref{Dix, cartan}, for $\nu \in (-2s,2s)$ and $\lambda, \mu \in (s^{-1/2},s^{1/2})$, becomes
	\begin{align*}
	\tau \left(E_{(0,s)} \left( B\right) \right)&\approx\int_{-2s}^{2s} \int_{s^{-1/2}}^{s^{1/2}} \int_{s^{-1/2}}^{s^{1/2}}  \tau \left(E_{(0,s)}\left[ \pi_{\lambda,\mu, \nu}(B) \right]  \right)    d\lambda\, d\mu\,d\nu\\
	&=\int_{-2s}^{2s} \int_{s^{-1/2}}^{s^{1/2}} \int_{s^{-1/2}}^{s^{1/2}} N_{\mathcal{B}_5}(s) d\lambda\, d\mu\,d\nu\\
	&\lesssim 2\sqrt{2} C \cdot \int_{-2s}^{2s} \int_{s^{-1/2}}^{s^{1/2}} \int_{s^{-1/2}}^{s^{1/2}} s^{5/2}  d\lambda\, d\mu\,d\nu \lesssim C^{'} \cdot s^{9/2}\,,\quad s \rightarrow \infty\,.
	\end{align*}
	 Then, for $\phi$ as in the statement, an application of Theorem \ref{Theorem 1.1,AR}, yields
	\begin{equation}
	\label{phi,ineq,cartan}
	\|\phi\left( B \right)\|_{L^{p}(\mathcal{B}_5) \rightarrow L^{q}(\mathcal{B}_5)} \lesssim C^{'} \cdot \sup_{s>0} \phi(s) \cdot  \tau \left( E_{(0,s)}(B)\right)^{\frac{1}{p}-\frac{1}{q}} \lesssim C^{'} \cdot  s^{\frac{9}{2r}}\,,
	\end{equation}
	where $\frac{1}{r}=\frac{1}{p}-\frac{1}{q}$, given that $1<p \leq 2 \leq q < \infty$. Finally, plugging $\phi$ as in \eqref{def_of_phi,eng} into \eqref{phi,ineq,cartan}, we get the Sobolev-type estimates  \eqref{Sob,cartan}.
\end{proof}
We note that, in the above example, the operator $B$ can be generalised as being any non-Rockland operator of the form 
\[
B=-(X_{1}^{2}+X_{2}^{2n_1}\pm X_{3}^{2n_2} \pm X_{4}^{2n_3} \pm X_{5}^{2n_4} \pm X_{4}^{-2n_5} \pm X_{5}^{-2n_6})\,,
\]
where $n_i \geq 2$, $n_i \in \mathbb{N}$ for $i=1,\ldots,6$, if we choose the positive sign of the corresponding left-invariant operator if $n_i$ is of the form $n_i=1+2k_i$, $k_i \geq 1$, and the negative sign, otherwise. 
\begin{rem}
\label{Sob,emb2,car}
For $b=0$ in the definition of the function $\phi$ in \eqref{def_of_phi,eng}, we obtain Sobolev-type embedding inequalities of the form
\[
\|f\|_{L^{q}(\mathcal{B}_5)} \leq C \| (I+B)^a f\|_{L^{p}(\mathcal{B}_5)}\,,
\]
for $a \geq \frac{9}{2r}$ and $1<p \leq 2 \leq q < \infty$.
\end{rem}
\begin{rem}(The $B$-heat equation)
\label{heat.car} For the operator $B$ as in Example \ref{multipliers,cartan} consider the heat equation given by
\begin{equation}
\label{heat.eq.car}
\partial_t u +Bu=0\,,\quad u(0)=u_0\,.
\end{equation}
Similarly to what has been done in Remark \ref{heat,eng}, we consider for $t>0$ the function $u(t)=u(t,x):=e^{-tB}u_0$ to be the solution of the initial value problem \eqref{heat.eq.car}, as well as the function $\phi=\phi(s):=e^{-ts}$ satisfying the assumptions given in Theorem \ref{Theorem 1.1,AR}. Reasoning as we did before in Remark \ref{heat,eng}, an application of Example \ref{multipliers,cartan} provides us with an estimate for time delay of the hear kernel of \eqref{heat.eq.car}, i.e., we have
\[
\|u(t, \cdot) \|_{L^q(\mathcal{B}_5)} \leq C^{'}_{p,q}t^{-\frac{9}{2}\left( \frac{1}{p}-\frac{1}{q}\right)}\|u_0\|_{L^p(\mathcal{B}_5)},\; 1<p \leq 2 \leq q< \infty\,.
\]
\end{rem}

\section*{Acknowledgement}
I would like to thank Professor Michael Ruzhansky for introducing me to this topic and for comments leading to improvements of the current work.


\begin{thebibliography}{CDR18}

\bibitem{AR20}
R.~Akylzhanov, M.~Ruzhansky.
\newblock $L^p$-$L^q$ multipliers on locally compact groups.
\newblock {\em J. Funct. Anal.}, 278 (2020), no. 3, 108324.

\bibitem{AK19}
F.~Almaki, V.~V. Kisil.
\newblock Geometric dynamics of a Harmonic oscillator, arbitary minimal uncertainly states and the smallest step 3 Nilponet {L}ie group.
\newblock{\em J. Phys. A: Math. Theor.}, 52 (2019), 025301.

%\bibitem[BFKG12]{BFKG12}
%H.~Bahouri, C.~Fermanian-Kammerer and I.~Gallagher.
%\newblock{\em Phase-space analysis and pseudo-differential calculus on the Heisenberg group},
%\newblock Ast\'{e}risque, \textbf{342}, 2012. See also revised version of March 2013 of arXiv:0904.4746.

\bibitem{BBR96}
P.~Boggiatto, E.~Buzano and L.~Rodino.
\newblock {\em Global hypoellipticity and spectral theory}, volume 92 of {\em Mathematical Research}.
\newblock Akademie Verlag, Berlin, 1996. 

\bibitem{BGR10}
U.~Boscain, J.~P. Gauthier and F.~Rossi.
\newblock Hypoelliptic heat kernel over 3-step nilpotent Lie groups.
\newblock{ \em J. Math. Sci.}, Vol. 199, No. 6, 2014.

\bibitem{CDR21}
D. Cardona, J. Delgado and M. Ruzhansky.
$L^p$-bounds for pseudo-differential operators on graded Lie groups. to appear in {\em J. Geom. Anal.}	arXiv:1911.03397

\bibitem{Cha20}
M.~Chatzakou.
\newblock On $(\lambda,\mu)$-classes on the Engel group. Advances in Harmonic Analysis and Partial Differential Equations, Trends Math., Birkh\"auser/Springer, (2020), 37-49. Cham. 

\bibitem{Cha21}
M.~Chatzakou.
\newblock Quantizations on the Engel and the Cartan groups. {\em J. Lie Theory}, 31 (2021), 517-542. 

\bibitem{CDR18}
M.~Chatzakou, J.~Delgado and M.~Ruzhansky.
\newblock On a class of anharmonic oscillators.
\newblock{ \em arXiv:2485649}, 2018. to appear in {\em  J. Math. Pures Appl.}

%\bibitem[CGGP92]{CGGP92}
%M.~Christ, D.~Geller, P.~Glowacki and L.~Polin.
%\newblock { Pseudodiffererential operators on groups with dilations}.
%\newblock {\em Duke Math. J.}, \textbf{68}, 379-423, 1998.

\bibitem{CG90}
L.~J. Corwin, F.~P. Greenleaf.
\newblock Representations on nilpotent Lie groups and their applications.
\newblock {\em Part I}, volume 18 of {\em Cambridge Studies in Advanced Mathematics}. Cambridge University Press, Cambridge, 1990. Basic theory and examples.

\bibitem{DR19}
J. Delgado, M. Ruzhansky.
$L^p$-bounds for pseudo-differential operators on compact Lie groups. 
{\em J. Inst. Math. Jussieu}, 18 (2019), 531–559.
 
% \bibitem[Dix81]{Dix81}
% J.~Dixmier.
% \newblock {\em Von Neumann algebras}. Amsterdam ; New York : North-Holland Pub. Co, 1981. 
 
\bibitem{Dix57}
J.~Dixmier.
\newblock Sur les repr\'{e}sentations unitaires des groupes de Lie nilpotents, volume III of {\em Canad. J. Math.}, 10, 321-348, 1957.

%\bibitem[Dix77]{Dix77}
%J.~Dixmier.
%\newblock {\em $C^{*}$ algebras}.
%\newblock North-Holland Publishing Co., Amsterdam, 1977.
%\newblock Translated from the French by Francis Jellet, North-Holland %Mathematical library, Vol.15.

\bibitem{Dix81}
J.~Dixmier.
\newblock {\em Von Neumann algebras}, volume 27 of North-Holland Mathematical library. North-Holland Publishing Co., Amsterdam, 1981.

%\bibitem[Dyn76]{Dyn76}
%A.~S. Dynin.
%\newblock An algebra of pseudodifferential operators on the Heisenberg groups. Symbolic calculus,
%\newblock {\em Dokl. Akad. Nauk SSSR}, \textbf{227}, 792-795, 1976.

%\bibitem[FR14]{FR14}
%V.~Fischer and M.~Ruzhansky.
%\newblock A pseudo-differential calculus on graded nilpotent Lie groups.
%\newblock{\em In Fourier analysis},Trends Math., pages 107-132. 
%\newblock Birkh\"auser/Springer, Cham, 2014.

\bibitem{FR16}
V.~Fischer, M.~Ruzhansky.
\newblock {\em Quantization on nilpotent {L}ie groups}, volume 314 of {\em
  Progress in Mathematics}.
\newblock Birkh\"auser/Springer, [Open access book], 2016.

\bibitem{FR17}
V.~Fischer, M.~Ruzhansky.
Sobolev spaces on graded groups.
{\em Ann. Inst. Fourier (Grenoble)}, 67 (2017), 1671-1723.

\bibitem{FR21}
V.~Fischer, M.~Ruzhansky.
Fourier multipliers on graded Lie groups.
{\em Colloq. Math.}, 165 (2021), 1-30.

%\bibitem[Fol94]{Fol94}
%G.B.~Folland.
%\newblock Meta-Heisenberg groups, Fourier analysis (Orono, ME, 1992), Lecture notes in Pure and Appl. Math.,
%\newblock \textbf{157}, 121-147, Dekker, New York, 1994.

\bibitem{FS82}
G. B.~ Folland, E.~Stein.
\newblock {\em Hardy spaces on Homogeneous groups}.
\newblock Mathematical Notes, {28}.
Princeton University Press, 1982.

\bibitem{Fol75}
G. B.~ Folland.
\newblock Subelliptic estimates and function spaces on nilpotent {L}ie groups.
\newblock {\em Ark. Mat.}, 13(2):161--207, 1975.

\bibitem{FS74}
G. B.~Folland, E. M.~Stein.
\newblock {Estimates for the $\overline{\partial}_{b}$ complex and analysis on the Heisenberg group.}
\newblock{\em Comm. Pure Appl. Math.}, 27, 1974, 429-522.

\bibitem{HK16}
A.~Hassannezhad, G.~Kokarev.
\newblock Sub-Laplacian eigenvalue bounds on sub-Riemannian manifolds. \newblock {\em Ann. Scuola Norm. Sup Pisa Cl. Sci} ,XVI(4):1049-1092, 2016.

%\bibitem[H\"{o}r85]{Hor85}
%L.~H\"{o}rmander.
%\newblock {\em The Analysis of linear partial differential operators}, vol. III,  \newblock {\em Springer-Verlang}, 1985. 

\bibitem{Hor60}
L.~H\"{o}rmander.
\newblock {Estimates for translation invariant operators in $L^{p}$ spaces}.
\newblock {\em Acta Math.}, 104:93–140, 1960.

\bibitem{RS77}
L. P.~Rothschild, E. M.~ Stein.
\newblock{Hypoelliptic differential operators and nilpotent groups.}
\newblock {\em Acta Math.}, 137, 1976, 3-4, 247-320.


\bibitem{RR18}
D.~Rottensteiner, M.~Ruzhansky.
\newblock {Harmonic and anharmonic oscillators on the Heisenberg group.}
\newblock {\em arXiv:1812.09620v1}, 2018.

%\bibitem[RT10]{RT10}
%M.~Ruzhansky and V.~Turunen.
%\newblock{\em Pseudo-differential operators and symmetries: Background analysis and advanced topics.}
%\newblock{Pseudo-Differential Operators: Theory and Applications},\textbf{2}, Birkh\"{a}user, Verlag, 2010.


%\bibitem[Shu01]{Shu01}
%M.A.~ Shubin.
%\newblock {\em Pseudodifferential operators and spectral theory. }
%\newblock Springer-Verlag, Berlin, second edition, 2001. 
%\newblock Translated from the 1978 Russian original by Stig I. Anderson.

%\bibitem[Tayl84]{Tayl84}
%M.E.~ Taylor.
%\newblock {\em Noncommutative microlocal analysis. I. Mem.}
%\newblock Amer. Math. Soc. \textbf{52}, 1984.
 
%\bibitem{Tit58}
%E.C.~ Titchmarsh 
%\newblock {\em Eigenfunction expansions associated with second-order differential equations}. Oxford University Press, Oxford, 1958.



\end{thebibliography}
 \end{document}